\documentclass[12pt]{amsart}
\usepackage{amssymb}

\begin{document}

\newcommand{\wt}{\widetilde}
\newcommand{\Aut}{\mbox{{\rm Aut}$\,$}}
\newcommand{\ul}{\underline}
\newcommand{\ol}{\overline}
\newcommand{\lr}{\longrightarrow}
\newcommand{\bc}{{\mathbb C}}
\newcommand{\bp}{{\mathbb P}}
\newcommand{\cf}{{\mathcal F}}
\newcommand{\ce}{{\mathcal E}}
\newcommand{\co}{{\mathcal O}}
\newcommand{\cg}{{\mathcal G}}
\newcommand{\hra}{\hookrightarrow}

\newtheorem{theorem}{Theorem}[section]
\newtheorem{lemma}[theorem]{Lemma}
\newtheorem{proposition}[theorem]{Proposition}
\newtheorem{corollary}[theorem]{Corollary}
\newtheorem{definition}[theorem]{Definition}
\newtheorem{remark}[theorem]{Remark}

\title{Universal families on moduli spaces of principal
bundles on curves}

\author{V. Balaji}

\address{Chennai Mathematical Institute, 92, G.N. Chetty Road,
Chennai 600017, India}

\email{balaji@cmi.ac.in}

\author{I. Biswas}

\address{School of Mathematics, Tata Institute of Fundamental
Research, Homi Bhabha Road, Mumbai 400005, India}

\email{indranil@math.tifr.res.in}

\author{D. S. Nagaraj}

\address{The Institute of Mathematical Sciences, CIT
Campus, Taramani, Chennai 600113, India}

\email{dsn@imsc.res.in}

\author{P. E. Newstead}

\address{Department of Mathematical Sciences, The University of
Liverpool, Peach Street, Liverpool, L69 7ZL, England}

\email{newstead@liverpool.ac.uk}

\date{}

\thanks{All authors are members of the international research group 
VBAC (Vector Bundles on Algebraic Curves), which is partially supported by 
EAGER (EC FP5 Contract no. HPRN-CT-2000-00099) and by EDGE (EC FP5 
Contract no. HPRN-CT-2000-00101). The fourth author would also like to
thank the Royal Society and CSIC, Madrid, for supporting a visit to CSIC 
during which his contribution to this work was completed.}

\begin{abstract}
Let $H$ be a connected semisimple linear algebraic group
defined over $\mathbb C$ and
$X$ a compact connected Riemann surface of genus at least three.
Let ${\mathcal M}'_X(H)$ be the moduli space
parametrising all topologically trivial stable
principal $H$-bundles over $X$ whose automorphism group
coincides with the centre of $H$. It is a Zariski open dense subset
of the moduli space of stable principal $H$-bundles. We prove that
there is a universal principal
$H$-bundle over $X\times {\mathcal M}'_X(H)$ if and only
if $H$ is an adjoint group (that is, the centre of $H$ is trivial).
\end{abstract}

\maketitle

\section{Introduction}

Let $X$ be a compact connected Riemann surface of genus at least
three. Let ${\mathcal M}_X(n,d)$ denote the moduli space of all
stable vector bundles over $X$ of rank $n$ and degree $d$, which is a
smooth irreducible quasiprojective variety defined over $\mathbb C$.
A vector bundle $\mathcal
E$ over $ X\times {\mathcal M}_X(n,d)$ is called
\textit{universal} if for every point
$m\in {\mathcal M}_X(n,d)$, the restriction of $\mathcal E$ to
$X\times \{m\}$ is in the isomorphism class of
holomorphic vector bundles over $X$ defined by $m$. A
well-known theorem says that there is a universal vector
bundle over $X\times {\mathcal M}_X(n,d)$ if and only if $d$ is
coprime to $n$ (see \cite{Ty} for existence in the coprime case,
\cite{Ra} for non-existence in the non-coprime
case and \cite{Ne} for a topological version of non-existence 
in the case $d=0$).

Let $H$ be a connected semisimple linear algebraic group defined
over the field of complex numbers.
Ramanathan extended the notion of (semi)stability to principal
$H$-bundles and constructed moduli spaces for stable
principal $H$-bundles over $X$
\cite{Ra1, Ra2, Ra3}. The construction works for any given topological
type, yielding a moduli space which is an irreducible quasiprojective
variety defined over $\mathbb C$. We are concerned here with the case
of topologically trivial stable principal $H$-bundles.
Let ${\mathcal M}_X(H)$ denote the moduli space of
topologically trivial stable principal $H$-bundles over $X$.

Let $Z(H)\, \subset\, H$ be the centre. For any $H$-bundle $E_H$, the
group $Z(H)$ is contained in the automorphism group $\text{Aut}(E_H)$.
Let
$$
{\mathcal M}'_X(H)\, \subset \, {\mathcal M}_X(H)
$$
be the subvariety consisting of all $H$-bundles $E_H$ over $X$ with the
property $\text{Aut}(E_H)\, \cong\, Z(H)$. It is known that
${\mathcal M}'_X(H)$ is a dense Zariski-open subset contained in
the smooth locus of ${\mathcal M}_X(H)$.

A principal $H$-bundle $E$ over $X\times {\mathcal M}'_X(H)$
will be called a
\textit{universal bundle} if for every
point $m\in {\mathcal M}_X(n,d)$, the
restriction of $\mathcal E$ to $X\times \{m\}$ is in the isomorphism
class of stable
$H$-bundles over $X$ defined by the point $m$ of the moduli space.

The following theorem is the main result proved here:

\begin{theorem}\label{theorem0}
There is a universal $H$-bundle over $X\times
{\mathcal M}'_X(H)$ if and only if $Z(H) \, =\, e$.
\end{theorem}

\begin{remark}\label{sl}
{\rm Non-existence for $H={\rm SL}(n,{\mathbb C})$
was previously known
\cite{Ra, Ne}, as also was existence for $H={\rm PGL}(n,{\mathbb C})$.}
\end{remark}

{\it Acknowledgments}. We thank A. King, A. Nair and P. Sankaran for
useful discussions.

\section{Existence of universal bundle}
We begin this section by recalling very briefly
certain facts from \cite{basa}. Unless otherwise specified,
all bundles and sections considered will be algebraic.

\begin{remark}\label{hb}
{\rm Let $E$ be a principal $G$-bundle over $X$, where $G$
is a reductive linear algebraic group
defined over $\mathbb C$, and $H \subset G$ is a Zariski closed
semisimple
subgroup. For any variety $Y$ equipped with an action of $G$, the
fibre bundle $(E\times Y)/G$
over $X$ associated to $E$ will be denoted by $E(Y)$.
{\renewcommand{\labelenumi}{{\rm (\alph{enumi})}}}
\begin{enumerate}
\item There is a natural action of the group $\Aut_GE$, defined
by all
automorphisms of $E$ over the identity map of $X$ that
commute with the action of $G$, on $\Gamma (X, E(G/H))$ (the
space of all holomorphic sections of the fibre bundle
$E(G/H) \,=\, E/H$ over $X$)
and the orbits correspond to the equivalence classes of
$H$-reductions of $E$ with two reductions
being equivalent if the corresponding principal $H$-bundles are
isomorphic.

\item Let $G = {\rm GL}(n,\bc )$, and let $\phi: H \hookrightarrow G$
be a faithful
representation of the semisimple group $H$. 
Let $Q$ denote the open subset of semistable principal $G$-bundles
(or equivalently, of trivial determinant semistable
vector bundles of rank $n$) of the usual ``Quot scheme'',
and
let $Q(\phi)$ be the ``Quot scheme'' which parametrises pairs
of the form $(E',s)$,
where $E'$ is a principal $G$-bundle and $s$ is a reduction of structure
group of $E'$
to $H$. Then $Q(\phi)$ is in a sense a ``relative Quot
scheme''. As is clear from the definition and the notation, this
scheme is dependent on the choice of the inclusion $\phi: H \hra
G$ (for details see \cite{basa, Ra2, Ra3}).

One also has a {\it tautological} sheaf on $X \times Q$ which in fact
is a vector bundle. We denote by $\ce$ the {\it associated
tautological} principal $G$-bundle on $X \times Q$.

Recall that the moduli space of principal $G$-bundles ($G =
{\rm GL}(n, {\mathbb C})$) is realised
as a good quotient of $Q$ by the action of a reductive group
$\cg$. We may also assume that the group $\cg$ is with trivial
centre (see for example \cite{drezet}).
\end{enumerate}}
\end{remark}

\begin{remark}
{\rm It is immediate that the action of $\cg$ on $Q$ lifts to an
action on $Q(\phi)$, where $\cg$, $Q$ and $Q(\phi)$ are defined
in Remark \ref{hb}.

We have a morphism
\begin{equation}\label{d.psi}
\psi \,:\, Q(\phi) \, \lr\, Q\, ,
\end{equation}
which sends any $H$-bundle $E'$ to the ${\rm GL}(n,\bc)$-bundle
obtained by extending the structure group of $E'$ using the
homomorphism $\phi$. In fact, $\psi$ is a $\cg$-equivariant
{\it affine morphism} (see \cite{basa}).}
\end{remark}

Continuing with the notation in the above two remarks,
consider the $\cg$-action on $Q(\phi)$ (defined in
Remark \ref{hb}(2)) with the linearisation
induced by the affine $\cg$-morphism $\psi$ in \eqref{d.psi}.
Since a good quotient of
$Q$ by $\cg$ exists and since $\psi$ is an affine
$\cg$-equivariant map, a good quotient $Q(\phi)// \cg$ exists
(see \cite[Lemma 5.1]{Ra3}).

Moreover by the universal property of categorical quotients,
the canonical morphism
\begin{equation}\label{opsi}
\ol{\psi} \,: \,Q(\phi)//\cg \,\lr\, Q//\cg
\end{equation}
given by $\psi$ is also {\it affine}.

\begin{theorem}\label{moduli}
Let ${\mathcal M}_X(H)$ denote the scheme
$Q(\phi)//\cg$ (see \eqref{opsi}). Then this scheme is the coarse
moduli scheme of semistable $H$-bundles. Further, the scheme
${\mathcal M}_X(H)$ is projective, and
if $H \hookrightarrow {\rm GL}(V)$ is a faithful representation,
then the canonical morphism
$$\ol{\psi}\, :\, {\mathcal M}_X(H) \,\lr\, {\mathcal
M}_X({\rm GL}(V))\, =\, Q//\cg
$$
is finite. 
\end{theorem}

Let $Q(\phi)^{s}$ be the open subscheme of $ Q(\phi)$ consisting of
stable $H$-bundles.

\begin{lemma}\label{free}
Let $Q_H' \subset Q(\phi)^s $ be the subset
parametrising all stable $H$-bundles
whose automorphism group is $Z(H)$. Then the action
of $\cg$ on the subset $Q_H'$ is free, and
furthermore, the quotient morphism $Q_H' \lr
Q_H'/{\cg}$ is a principal $\cg$-bundle. In fact, $Q_H'/{\cg}$ is
precisely the Zariski open subset ${\mathcal M}'_X(H)$ (the variety
${\mathcal M}'_X(H)$ is defined in the introduction).
\end{lemma}

\begin{proof}
For any point $E_H \in Q(\phi)$, the isotropy
subgroup of $E_H$ for the action of $\cg$ coincides with
${\rm Aut}(E_H)/Z(H)$.
This can be seen as follows: firstly, the point $E_H$ is a pair
$(E,s)$, where $E \in Q$ and $s$ is a reduction of
structure group to $H$ of the $G$-bundle $E$. It
is well known for the action of $\cg$ on $Q$ that the isotropy
at $E$ is precisely the group ${\rm Aut}(E)/Z(G)$. From this
it is easy to see that the isotropy of $E_H$ for the action of
$\cg$ on $Q(\phi)$ is the group
${\rm Aut}(E,s)/Z(H) = {\rm Aut}(E_H)/Z(H)$.

Hence it follows that the action of $\cg$ on the open subset $Q_H'$
is free, and the proof of the lemma is complete.
\end{proof}

We remark that, at least when $H$ is of adjoint
type, $Q_H'$ is a non-empty open subset of $Q(\phi)^{s}$. Openness is
easy and can be seen for example from \cite[Theorem II.6
(ii)]{faltings}. Non-emptiness follows from Proposition \ref{codim}
below.

We prove a proposition on semisimple groups, possibly known to
experts but which we could not locate in any standard text.
 
\begin{proposition}\label{faithrep}
Let $H$ be a semisimple algebraic group. Then
$H$ has a faithful, irreducible representation $\phi:H \lr {\rm GL}(V)$
if and only if the centre of $H$ is cyclic.
\end{proposition}

\begin{proof}
One way the implication is easy, namely suppose
that a faithful
irreducible representation exists, then the centre $Z(H)$
of $H$ is cyclic. 
To see this, first note
that under the representation $\phi$, the centre
$Z(H)$ maps to a subgroup which commutes
with all elements of $\phi(H)$. Since $\phi$ is irreducible, this
implies by Schur's Lemma that $\phi(Z(H)) \subset
Z({\rm GL}(V))$, where $Z({\rm GL}(V))$ is the centre
of ${\rm GL}(V)$. Observe
further that since $H$ is semisimple, we have $\phi(H) \subset {\rm
SL}(V)$. Hence, $Z(H) \subset Z({\rm SL}(V))$ and is therefore cyclic.

The contrapositive statement is harder to prove. We proceed as
follows: Let $H'$ be the simply connected cover of $H$, and let
$\overline H$ be the associated adjoint group, namely $H/Z(H)$. Let
$\Lambda$ (respectively, $\Lambda_R$) be the weight lattice
(respectively, root
lattice). In other words, by the Borel-Weil theorem $\Lambda =
{\mathcal X}(B)$ and $\Lambda_R = {\mathcal X}(\overline B)$, where 
$B$, $\overline B$ are fixed Borel subgroups of $H$, $\overline H$.
Then, from the exact sequence
\[
e \lr Z(H) \lr H \lr \overline H \lr e
\]
we see that $\Lambda/\Lambda_R \simeq {\mathcal X}(Z(H))$. In other words,
the quotient group $\Lambda/\Lambda_R$ is cyclic of order $m$.

Let $\overline \lambda$ be a generator of the cyclic
group $\Lambda/\Lambda_R$.
Then, by an action of the Weyl group we may assume that the coset
representative $\lambda \in \Lambda$ is actually a dominant weight.

Suppose that the root lattice has the following decomposition
(corresponding to the simple components of $\overline H$):
\[
\Lambda_R = \bigoplus_{i=1}^\ell \Lambda^{i}\, .
\]
Then, by possibly adding dominant weights from the $\Lambda^{i}$,
we may assume that $(m \cdot \lambda) \in \Lambda_R$ has
all its direct sum components $\lambda_i \neq 0$,
$i\, \in\, [1\, ,\ell]$, where $\lambda$ is as above.

For this choice of $\lambda \in \Lambda$ let $V_{\lambda}$ be the
corresponding $H$-module given by
\[
\phi_{\lambda} : H \lr {\rm GL}(V_{\lambda})\, .
\]
Then one knows that $\phi_{\lambda}$
is an irreducible representation of $H$.

We claim that the representation
$\phi_{\lambda}$ is even {\it faithful}. Suppose that this is
not the case.
Let $K_\lambda \,:= \,{\rm kernel}(\phi_{\lambda}) \,\neq\, e$.

Firstly, $K_\lambda
\subset Z(H)$. To see this, let $K'$ be the inverse
image of $K_\lambda$ in the simply connected cover $H'$ of $H$.
Then the choice of
$\lambda$ so made that its simple components are non-zero in fact
forces the following: Suppose that $H' = H_1 \times \cdots
\times H_\ell$ is the decomposition of $H'$ into its almost simple
factors. (We recall that a semisimple algebraic group is called
almost simple if the quotient of it by its centre is simple.)
Then the normal subgroup $K'$ in its decomposition in $H'$ is
such that $K_i\, :=\, K'\bigcap H_i$
are proper normal subgroups of $H_i$. In particular,
$K_i \subset Z(H_i)$ for all $i\, \in\, [1\, ,\ell]$. This implies
that $K' \subset Z(H')$ and hence, $K_\lambda \subset Z(H)$.

Note that the dominant character $\lambda$ is non-trivial on the
generator of the centre $Z(H)$ because ${\mathcal X}(Z(H)) 
= \Lambda/\Lambda_R$. Now $K_\lambda \subset Z(H)$ and $Z(H)$
cyclic implies that $\lambda$ is non-trivial on
the generator of $K_\lambda$ as well. This contradicts the
fact that $K_\lambda = {\rm kernel}(\phi_{\lambda})$. This proves
the claim. Therefore, the proof of the proposition is complete.
\end{proof}

If $H$ is of adjoint type, then by
Proposition \ref{faithrep}
we can choose the inclusion $\phi: H \hra G$ in
Remark \ref{hb}(2) to be an irreducible representation.
Henceforth, $\phi$ will be assumed to be irreducible.

Let $E$ be a stable $H$-bundle of trivial topological
type. Recall that one can realise $E$ from a
unique, up to an inner conjugation, irreducible
representation of $\pi_1(X)$ in a
maximal compact subgroup of $H$ (see \cite{Ra1}).
For notational convenience, we will always
suppress the base point in the notation of fundamental group.
Denote by $M(E)$ the
Zariski closure of the image of $\pi_1(X)$ in $H$.

\begin{proposition}\label{codim}
Let $H$ be of adjoint type and let $\phi$ be a faithful
irreducible representation of $H$ in $V$ (see
Proposition \ref{faithrep}). Let $U \subset
Q(\phi)^{s}$ be the subset defined as follows:
\[
U = \{E \in Q(\phi)^{s}\, |\, M(E) = H\}.
\]
Then $U$ is non-empty and is contained in the subset $Q_H'$ of
stable $H$-bundles which have trivial automorphism group. (Since
$H$ is of adjoint type, $Z(H)=e$.)
\end{proposition}

\begin{proof}
Since $\phi$ is irreducible, using
Lemma 2.1 of \cite{Su} we conclude that there is an irreducible
representation
$$
\rho\, :\, \pi_1(X)\,\lr\, H
$$
which have the property that the composition
$\phi \circ \rho :\pi_1(X)\lr G := \text{GL}(V)$
continues to remain irreducible. This
implies that $M(E_{\rho}) = H$ by the
construction in \cite[Lemma 2.1]{Su} and hence $E_{\rho} \in U$, i.e.,
$U$ is non-empty.

Let $E \in U$. Then by the definition of $U$, there
exists a representation $\rho$ of $\pi_1(X)$ in $H$ such that $E
\simeq E_{\rho}$, where $E_{\rho}$ is the flat principal
$H$-bundle given
by $\rho$. Observe that $E_{\rho}$ is a stable $H$-bundle and
the associated $G$-bundle is also stable. Hence all
the automorphisms of this associated $G$-bundle
lie in $Z(G)$. Since $H$ is of adjoint type, it follows
that the $H$-bundle $E_{\rho}$ has no non-trivial
automorphisms. Hence it follows that $U \subset Q_H'$.
\end{proof}

We have the following theorem on existence of universal families.

\begin{theorem}\label{family}
Let $H$ be a group of adjoint type and ${\mathcal
M}'_X(H)$ the Zariski open subset of ${\mathcal
M}_X(H)$ defined in the introduction. Then there exists a
universal family of principal $H$-bundles on $X \times {\mathcal
M}'_X(H)$.
\end{theorem}

\begin{proof}
We first observe that the variety ${\mathcal
M}'_X(H)$ is precisely the image of $Q_H'$
under the quotient map for the action of $\mathcal G$.
We recall that $Q_H'$ is nonempty by Proposition \ref{codim}.

Since $Z(H) = e$, it follows that $H \subset {\ol G} \,:=\,
G/Z(G)$, where $G ={\rm GL}(V)$ is as in Remark \ref{hb}(2).

Consider the tautological $G$-bundle $\ce$ on $X \times Q^s$, and let
$\ol {\ce}$ be the corresponding $\ol G$-bundle
obtained by extending the structure group. Then it is well-known
that the adjoint universal bundle $\ol {\ce}$ descends to the quotient
${\mathcal M}_X(G)^s$ (see \cite{drezet}).

We follow the same strategy for ${\mathcal M}'_X(H)$ as well.
Consider the pulled back ${\ol G}$-bundle
$(\text{Id}_X \times {\psi})^{*} \ol {\ce}$
on $X\times Q_H'$, where
$\psi$ is the map in \eqref{d.psi} and ${\ol G} \,:=\,
G/Z(G)$.
The action of $\cg$ on $Q_H'$ is free by Lemma \ref{free}. Therefore,
the quotient $Q_H' \lr {\mathcal M}'_X(H)$ is a principal
${\cg}$-bundle. Further, the action of $\cg$ lifts to the
tautological bundle $(\text{Id}_X \times {\psi})^{*} \ol {\ce}$. In
particular,
the principal $\ol G$-bundle $(\text{Id}_X\times {\psi})^{*}
\ol {\ce}$ descends to a principal $\ol G$-bundle
over $X \times {\mathcal M}'_X(H)$.

Let us denote this descended $\ol G$-bundle over
$X \times {\mathcal M}'_X(H)$ by $\ol {\ce}_0$.

Let $\pi \,:\,
(\text{Id}_X\times \psi)^*{\ce}(G/H)\,\lr\,
(\text{Id}_X\times \psi)^*{\ol {\ce}}({\ol G}/H)$ be the natural
map induced by the projection $G/H \lr {\ol G}/H$, where
$\psi$ is the map in \eqref{d.psi}. We note that the
universal $H$-bundle over $X \times U$, where $U$ is defined in
Proposition \ref{codim}, is a reduction of structure group of
the pulled back $G$-bundle $(\text{Id}_X\times \psi)^*\ce$. Let
$$
\sigma \,:\, X\times U\, \longrightarrow\,
(\text{Id}_X\times \psi)^*{\ce}(G/H)
$$
be the section giving this reduction of structure group. Then
the composition ${\pi}\circ \sigma$
is a section of $(\text{Id}_X\times \psi)^*{\ol {\ce}}({\ol G}/H)$
over $X\times U$.

Since $H$ is semisimple, a lemma of Chevalley says that there
is a ${\ol G}$-module $W$ and an element
$w\in W$ such that $H$ is precisely the
isotropy subgroup for $w$ (see \cite[p. 89, Theorem 5.1]{Bo}).
Therefore, ${\ol G}/H$ is identified with
the closed ${\ol G}$-orbit in $W$ defined by $w$.
Then we see that ${\ol {\ce}}({\ol G}/H) \hra
{\ol{\ce}}(W)$. We may therefore view the section
${\pi}\circ \sigma: X\times U\longrightarrow
(\text{Id}_X\times \psi)^*{\ol {\ce}}({\ol G}/H)$ as a section of
the vector bundle $(\text{Id}_X\times \psi)^*{\ol {\ce}}(W)$
over $X\times U$.

Since $(\text{Id}_X \times {\psi})^{*} \ol {\ce}$ descends to the
${\ol G}$-bundle ${\ol {\ce}}_{0}$ over $X \times {\mathcal M}'_X(H)$,
it follows that the associated vector bundle
$(\text{Id}_X \times {\psi})^{*} \ol {\ce}(W)$
also descends to $X \times
{\mathcal M}'_X(H)$. Clearly this vector bundle is nothing but the
associated vector
bundle ${\ol {\ce}}_{0}(W)$, associated to the $\ol G$-bundle ${\ol
{\ce}}_{0}$ on $X \times {\mathcal M}'_X(H)$ for the $\ol G$-module
$W$.
 
Since each point of ${\mathcal M}'_X(H)$ represents an isomorphism
class of stable $H$-bundle, it follows that set theoretically the
reduction section ${\pi}\circ \sigma \,\in \,
{\Gamma}((\text{Id}_X\times \psi)^*{\ol {\ce}}(W))$
descends to a section on $X \times {\mathcal M}'_X(H)$ of the
descended vector bundle ${\ol {\ce}}_{0}(W)$.

We now appeal to \cite[Proposition 4.1]{knr}, which implies that the
section ${\pi}\circ \sigma$ in fact descends to give a holomorphic
section
of ${\ol{\ce}}_{0}(W)$ over $X \times {\mathcal M}'_X(H)$.

Again set-theoretically, the image of this section
of ${\ol{\ce}}_{0}(W)$ lies in ${\ol
{\ce}}_{0}({\ol G}/H)\, \subset\, {\ol{\ce}}_{0}(W)$.
As before, from \cite[Proposition 4.1]{knr}
it follows that ${\pi}\circ \sigma$ gives
a reduction of structure group to $H$
of the descended ${\ol G}$-bundle
${\ol{\ce}}_{0}$. The $H$-bundle over $X \times {\mathcal M}'_X(H)$
obtained this way is the required universal $H$-bundle.
This completes the proof of the theorem.
\end{proof}

\section{Nonexistence of universal bundle}

Let $H$ be a complex semisimple linear algebraic group, and let
$K\, \subset\, H$ be a maximal compact subgroup. The Lie algebra
of $K$ will be denoted by $\mathfrak k$. A homomorphism
$$
\rho\, :\, \pi_1(X)\, \longrightarrow\, K
$$
is called \textit{irreducible} if no nonzero vector in $\mathfrak k$
is fixed by the adjoint action of the
subgroup $\rho (\pi_1(X))\, \subset \, K$ on $\mathfrak k$.
Let $\text{Hom}^{\text{irr}}(\pi_1(X), K)$ denote the
space of all irreducible homomorphisms such that the corresponding
$K$-bundle is topologically trivial. So any homomorphism
in $\text{Hom}^{\text{irr}}(\pi_1(X), K)$ is induced by a
homomorphism from $\pi_1(X)$ to the universal cover of $K$.

Let $g$ denote the genus of $X$. Assume that $g\, \geq\, 3$.
If we choose a basis
$$
\{a_1, \cdots ,a_g, b_1, \cdots ,b_g\} \,\subset\, \pi_1(X)
$$
such that 
\begin{equation}\label{pi-rep.}
\pi_1(X) \,=\, {\langle} a_1,\ldots ,a_g, b_1,\ldots ,b_g :
\prod_{i=1}^na_ib_ia_i^{-1}b_i^{-1} \rangle\, ,
\end{equation}
then $\text{Hom}^{\text{irr}}(\pi_1(X), K)$ gets identified
with a real analytic subspace of $K^{2g}$.

For any $\rho\, \in\, \text{Hom}^{\text{irr}}(\pi_1(X), K)$, the
principal $H$-bundle obtained by extending the structure group of the
principal $K$-bundle given by $\rho$
is stable. A theorem of Ramanathan says that
all topologically trivial stable $H$-bundles arise in this way, 
that is, the space of all
equivalence classes of irreducible homomorphisms of $\pi_1(X)$
to $K$
is in bijective correspondence with the space of all stable
$H$-bundles over $X$
\cite[Theorem 7.1]{Ra1}. More precisely (see the proof of
\cite[Theorem 7.1]{Ra1}), the real analytic space
underlying the moduli space ${\mathcal M}_X(H)$ parametrising
topologically trivial stable $H$-bundles is analytically
isomorphic to the quotient space
\[
R(\pi_1(X),H)\, :=\, \text{Hom}^{\text{irr}}(\pi_1(X), K)/K\, ,
\]
for the action constructed using
the conjugation action of $K$ on itself. Let
\begin{equation}\label{q-pr.-b00.}
q\, :\, \text{Hom}^{\text{irr}}(\pi_1(X), K) \,
\longrightarrow\, R(\pi_1(X),H) \,\cong\, {\mathcal M}_X(H)
\end{equation}
be the quotient map. The open
subset ${\mathcal M}'_X(H)\, \subset\, {\mathcal M}_X(H)$ is a smooth
submanifold and the restriction of the map in \eqref{q-pr.-b00.}
\begin{equation}\label{q-pr.-b.}
q\vert_{q^{-1}({\mathcal M}'_X(H))}\, :\,
q^{-1}({\mathcal M}'_X(H))\, \longrightarrow\, {\mathcal M}'_X(H)
\end{equation}
is a smooth principal $K/Z$-bundle,
where $Z$ is the centre of $K$. Note that $Z$ coincides
with the centre of $H$ (as $H$ is semisimple, its centre is
a finite group).

If ${\mathcal U}_H$ is a universal $H$-bundle over $X\times {\mathcal
M}'_X(H)$, then we have a reduction of structure group
of $\mathcal U_H$ to $K$ which is constructed using the correspondence
established in \cite{Ra1}
between stable $H$-bundles and irreducible flat $K$-bundles. Let
\begin{equation}\label{d.uk}
{\mathcal U}_K \,\lr\, X\times {\mathcal M}'_X(H)
\end{equation}
be the smooth principal $K$-bundle obtained from
${\mathcal U}_H$ this way. Consider the
principal $K/Z$-bundle over $X\times {\mathcal M}'_X(H)$,
where $Z\, \subset\, K$ is the centre, obtained
by extending the structure group of ${\mathcal U}_K$ using the
natural projection of $K$ to $K/Z$.
The restriction of this $K/Z$-bundle to
$p\times {\mathcal M}'_X (H)$, where $p$ is the base point in $X$
used for defining $\pi_1(X)$, is identified with the $K/Z$-bundle
in \eqref{q-pr.-b.}. To see this first note that the universal $K$-bundle
over $X\times {\rm Hom}(\pi_1(X),K)$ is obtained as a quotient
by the action of $\pi_1(X)$ on $\widetilde{X}\times
{\rm Hom}(\pi_1(X),K)\times K$, where $\widetilde{X}$ is the pointed
universal cover of $X$ for the base point $p$; the
action of $z\in \pi_1(X)$ sends any $(\alpha,\beta,\gamma)$
to $(\alpha ,z, \beta , \beta(z)^{-1}\gamma)$. From this
it follows that the restriction of this universal bundle
to $p\times {\rm Hom}(\pi_1(X),K)$ is canonically trivialized.
The above mentioned identification is constructed using
this trivialization.

Our aim is to show that no $K$-bundle over
${\mathcal M}'_X(H)$ exists that produces the $K/Z$-bundle
in \eqref{q-pr.-b.} by extension of structure group, provided
the centre $Z$ is non-trivial.

Let $F_3$ denote the free group on three generators. Fix
a surjective homomorphism
\begin{equation}\label{def.f}
f\, :\, \pi_1(X) \,\longrightarrow\, F_3
\end{equation}
that sends $a_i$, $1\leq i\leq 3$, to the $i$-th generator of $F_3$
and sends $a_i$, $4\leq i\leq g$, and $b_i$, $1\leq i\leq g$, to the
identity element, where $a_j$, $b_j$ are as in \eqref{pi-rep.}
(recall that $g\geq 3$).

Set
\begin{equation}\label{qmap}
R(F_3,H) \,:=\, \text{Hom}^{\text{irr}}(F_3, K)/K
\end{equation}
to be the equivalence classes of irreducible representations.
Note that $\text{Hom}(F_3, K) \simeq K^3$, and under this
identification the action
$$
\mu \,:\, \text{Hom}(F_3, K)\times K
\,\longrightarrow\, \text{Hom}(F_3, K)
$$
given by $\mu(\rho, A) = A^{-1}\rho A$ corresponds to the
simultaneous diagonal
conjugation action of $K$ on the three factors.

Since $F_3$ is a free group and $K$ is connected, any
homomorphism $\rho$ from $F_3$ to $K$ can be deformed to the trivial
homomorphism. This implies that the principal $H$-bundle corresponding
to the homomorphism
$$
\rho\circ f\, :\, \pi_1(X) \,\longrightarrow\, K
$$
is topologically trivial, where $f$ is defined
in \eqref{def.f}. Therefore, we have an embedding
\begin{equation}\label{i.map}
R(F_3,H) \,\longrightarrow\, R(\pi_1(X),H)\, \cong\,
{\mathcal M}_X(H)
\end{equation}
that sends any $\rho$ to $\rho\circ f$.
Let $R'(F_3,H)\, \subset\, R(F_3,H)$ be the inverse image of
the open subset ${\mathcal M}'_X(H)$ under the above map.
By Proposition \ref{codim} it follows that $R'(F_3,H)$ is a 
{\it non-empty} open subset of $R(F_3,H)$.

Fix a point $p \in X$, which will also be the base
point for the fundamental group.
Consider the restriction of the principal
$K$-bundle ${\mathcal U}_K$ in \eqref{d.uk}
to $p\times {\mathcal M}'_X(H)\, \hra\, X\times {\mathcal M}'_X(H)$
and denote it by ${\mathcal U}_{K,p}$. 

Let
\begin{equation}\label{1ga}
\gamma \,:\, p\times R'(F_3,H))\,\lr\, p\times{\mathcal M}'_X(H)
\end{equation}
be the map given by the embedding
$R'(F_3,H) \lr {\mathcal M}'_X(H)$
constructed in \eqref{i.map}.
Taking the pull back of the above defined
principal $K$-bundle ${\mathcal U}_{K,p}$ over
$p\times {\mathcal M}'_X(H)$ under the morphism $\gamma$
in \eqref{1ga}, 
\begin{equation}\label{eqdig}
\begin{matrix}
 \widetilde M && {\mathcal U}_{K,p} \\
 \Big\downarrow & & \Big\downarrow\\
p\times R'(F_3,H))
& \stackrel{\gamma}\lr & p\times {\mathcal M}'_X(H) \\
\end{matrix}
\end{equation}
we obtain a principal $K$-bundle $\widetilde{M} \, :=\,
{\gamma}^{*}{\mathcal U}_{K,p}$ on $R'(F_3,H)$ whose associated
$K/Z$-bundle (as before, $Z$ is the centre of $K$)
is precisely the $K/Z$-bundle
$q^{-1}( R'(F_3,H)) \,\longrightarrow\,
R'(F_3,H)$ constructed in \eqref{q-pr.-b.}.
In particular, $q^{-1}( R'(F_3,H))\, \simeq\, {\widetilde
M}/Z$.

Let $\text{Hom}^{\text{irr}}(F_3, K)'\, \subset\,
\text{Hom}^{\text{irr}}(F_3, K)$ be the inverse image of
$R'(F_3,H)$ under the quotient map in \eqref{qmap}. Let
\begin{equation}\label{map.q0}
q_0\, :\, \text{Hom}^{\text{irr}}(F_3, K)' \,\longrightarrow\,
R'(F_3,H)
\end{equation}
be the restriction of the quotient map in \eqref{qmap}.
So $q_0$ is the quotient for the conjugation action of $K$.
This map $q_0$ defines a principal $K/Z$-bundle
over $R'(F_3,H)$ which is evidently identified
with the $K/Z$-bundle ${\widetilde M}/Z\,\longrightarrow\,
R'(F_3,H)$ obtained from \eqref{eqdig}.

We will prove (by contradiction) that such a $K$-bundle
$\widetilde{M}$ does not exist. In other words, we will show
that there is no principal $K$-bundle over $R'(F_3,H)$ whose
extension of structure group is the $K/Z$-bundle in \eqref{map.q0}.

Assume the contrary, then we get the following diagram of
topological spaces and morphisms:
\begin{equation}\label{eqdig2}
\begin{matrix}
Z & \longrightarrow & K & \longrightarrow & K/Z \\
|| & {} & \Big\downarrow & & \Big\downarrow\\
 Z & \longrightarrow & \widetilde{M}
& \longrightarrow &q^{-1}_0( R'(F_3,H))\\
& & \Big\downarrow & & \Big\downarrow\\
& & R'(F_3,H) & = & R'(F_3,H)
\end{matrix}
\end{equation}

We will need a couple of results. The proof of the
non-existence of the bundle $\widetilde{M}$ will be
completed after establishing Proposition \ref{claim}.

\begin{lemma}\label{newstead}
There exist finitely many compact differentiable
manifolds $N_i$ and differentiable maps 
\[
f_i \,:\, N_i\, \lr\, K^3\, ,
\]
$1\,\leq\, i\,\leq\,d$, such that if $N \,:= \,
\bigcup_{i=1}^d f_i(N_i)$ then
\begin{enumerate}
\item the complement $K^3 \setminus N$ is contained
in $q^{-1}_0(R'(F_3,H))$, where $q_0$ is the projection
in \eqref{map.q0}.

\item Furthermore, $\dim K^3 - \dim N_i \,\geq\, 4$
for all $i\,\in\, [1\, ,d]$.
\end{enumerate}
\end{lemma}

\begin{proof}
Take any homomorphism $\rho\, :\, \pi_1(X)\,
\longrightarrow\, K$. Let $E_\rho$ denote the
corresponding polystable principal $H$-bundle over
$X$ \cite{Ra1}. The automorphism group of the
polystable $H$-bundle $E_\rho$ coincides with
the centraliser of $\rho(\pi_1(X))$ in $H$. Hence
if $\rho(\pi_1(X))$ is dense in $K$, then the
automorphism group of $E_\rho$ is the centre
$Z(H)\, \subset\, H$. Therefore, if
the topological closure $\overline{\rho(\pi_1(X))}$
of $\rho(\pi_1(X))$ is $K$,
and $E_\rho$ is topologically trivial, then
$E_\rho\, \in\, {\mathcal M}'_X(H)$.

Now assume that the centraliser
$C(\rho)\, \subset\, H$ of $\rho(\pi_1(X))$
in $H$ properly contains $Z(H)$. Since the
complexification of
$\overline{\rho(\pi_1(X))}$ is reductive, and the
centraliser of a reductive group is reductive,
we conclude that $C(\rho)$ is reductive.
Take a semisimple element $z\, \in\, (H\setminus Z(H))
\bigcap C(\rho)$ (since $C(\rho)$ is reductive
and larger than $Z(H)$ such
an element exists). Let $C_z\, \subset\, K$ be the
centraliser of $z$ in $K$. We have
$\overline{\rho(\pi_1(X))}\, \subset\, C_z$, and
$C_z$ is a proper subgroup of $K$ as $z\, \notin\,
Z(H)$. Also, $C_z$ contains a maximal torus of $K$.

Fix a maximal torus $T\, \subset\, K$. (Since any
two maximal tori are conjugate, any subgroup $C_z$ of the
above type would contain $T$ after an inner conjugation.)
Consider all proper Lie subgroups $M$ of $K$ satisfying the
following two conditions:
\begin{enumerate}
\item[(i)] $T\, \subseteq\, M$, and
\item[(ii)] there exists a semisimple
element $z\, \in\, H\setminus Z(H)$
such that $M$ is the centraliser of $z$ in $K$.
\end{enumerate}
Let $\mathcal S$ denote this collection of Lie subgroups of $K$.

The connected ones among $\mathcal S$ are
precisely the maximal compact subgroups of Levi
subgroups of proper
parabolic subgroups of $H$ containing $T$. Note that if $P$ is
a proper parabolic subgroup of $H$, then
$P\bigcap K$ is a maximal compact subgroup of a Levi
subgroup of $P$. All the connected ones among
$\mathcal S$ arise as $P\bigcap K$ for some proper parabolic
subgroup $P\, \subset\, H$ containing $T$.

This collection $\mathcal S$ is a finite set. To see this,
we first note that there are
only finitely many parabolic subgroups of $H$ that
contain $T$. For a proper parabolic subgroup $P\,\subset\,
H$ containing $T$, there are only finitely many Lie
subgroups of $K$ that have $P\bigcap K$ as the connected
component containing the identity element. Thus
$\mathcal S$ is a finite set.

Let $M_1, \cdots , M_d$ be the subgroups of $K$ that occur
in $\mathcal S$.

We will show that the codimension of each $M_i$ in $K$
is at least two. It suffices show that for a
maximal proper parabolic subgroup $P$ of $G$ containing $T$,
the codimension of $M\, :=\, P\bigcap K$ in $K$ is at least two.

To prove that the codimension of $M\, =\, P\bigcap K$
in $K$ is at least two,
let $\mathfrak k$ be the Lie algebra
of $K$, and let $\mathfrak h$ be the Lie algebra of $T$
(the Cartan subalgebra). Let
\[
{\mathfrak k} = {\mathfrak h} + \sum_{\alpha~{\rm a}~{\rm root}}
{\mathfrak k}^{\alpha} 
\]
be the root space decomposition of $\mathfrak k$. Let $\mathfrak m$ be a
Lie
algebra of $M$ which is a Lie subalgebra of $\mathfrak k$ containing
$\mathfrak h$. Then we have the decomposition
\[
{\mathfrak m} \,=\, {\mathfrak h} + \sum_{\alpha~{\rm a}~{\rm root}}
{\mathfrak m}^{\alpha} 
\]
for $\mathfrak m$ where each ${\mathfrak m}^{\alpha} $ is irreducible for
$\mathfrak h$, and hence has to coincide with one of the ${\mathfrak
k}^{\alpha}$. Then $\mbox{codim}_{_{\mathfrak k}}({\mathfrak m}) \geq 2$
since
if $\alpha$ is a root for the subalgebra $\mathfrak m$ then so is
$-{\alpha}$ (see \cite[p. 83, Corollary 4.15]{Ad}).

Consequently, the codimension of $M$ in $K$ is at least two.
Thus the codimension of each $M_i$, $i\, \in\, [1\, ,d]$, in $K$
is at least two.

For each $i\, \in\, [1\, ,d]$, let
\begin{equation}\label{dbfi}
\overline{f}_i \,:\, K \times M^3_i \,\lr\, K^3
\end{equation}
be the map
defined by $(x,(y_1, y_2,y_3)) \,\longmapsto\, (x y_1 {x}^{-1},
x y_w {x}^{-1}, x y_3 {x}^{-1})$. Consider the free action
of $M_i$ on $K \times M^3_i$ defined by
\[
z \cdot (x,(y_1,y_2,y_3)) = (x z^{-1}, (z y_1 z^{-1},
z y_2 z^{-1}, z y_3 z^{-1}))\, ,
\]
where $x\in K$ and $z,y_1,y_2, y_3\in M_i$.
The map $f_i$ in \eqref{dbfi} clearly factors through the quotient
\[
\frac {K \times {M_i}^t}{M_i}
\]
for the above action. Therefore, we have
\begin{equation}\label{dfi}
f_i \,:\, N_i\, :=\, \frac {K \times {M_i}^t}{M_i}\, \lr\, K^t
\end{equation}
induced by $f_i$.

To prove part (1) of the lemma, we recall
the earlier remark that for any homomorphism
$\rho' \,\in\, q^{-1}_0(R'(F_3,H))$ (the map
$q_0$ is defined in \eqref{map.q0}), the automorphism
group of the principal
$H$-bundle corresponding to $\rho'\circ f$
(the homomorphism $f$ is defined in \eqref{def.f})
coincides with $Z(H)$
if the image of $\rho'$ is dense in $K$. From the
properties of the collection $\{M_1,\cdots , M_d\}$
we conclude that
$$
N\, :=\, \bigcup_{i=1}^d f_i(N_i)\, \supseteq\,
q^{-1}_0(R'(F_3,H))^c\, ,
$$
where $f_i$ are defined in \eqref{dfi} and
$$
q^{-1}_0(R'(F_3,H))^c\, \subset\, K^3
$$
is the complement of $q^{-1}_0(R'(F_3,H))$ in $K^3$.

Therefore, proof of part (1) is complete. To prove
part (2), we note that
$$
\dim N_i \,=\, \dim K + 2\cdot \dim M_i
$$
for all $i\,\in\, [1\, ,d]$.
It was shown earlier that $\dim M_i \, \leq\, \dim K -2$.
Therefore,
$$
\dim K^3 - \dim N_i \, =\,
3\cdot \dim K - \dim N_i \, \geq\, 3\cdot \dim K
- 3\cdot \dim K +4 \, =\, 4\, .
$$
This completes the proof of the lemma.
\end{proof}

\begin{proposition}\label{claim}
Consider the $K/Z$-principal bundle $q^{-1}_0( R'(F_3,H)) \,
\longrightarrow\, R'(F_3,H)$, where $q_0$ is the projection in
\eqref{map.q0}. The induced homomorphism on
fundamental groups
$$
\pi_1(K/Z) \,\longrightarrow\,\pi_1(q^{-1}(R'(F_3,H)))
$$
obtained from the homotopy exact sequence is trivial.
\end{proposition}

\begin{proof}
Let $x_1, x_2, x_3 \in K$ be regular
elements; we recall that $x\in K$ is called regular
if the centralizer $C(x)= \{y \in K \,\vert\, yx=xy \}$ is a
maximal torus in $K$. Since the set of regular elements is dense in
$K$, and $q^{-1}_0(R'(F_3,H)) \,\subset\, K^3$ is a nonempty
open dense subset (this follows from Lemma \ref{newstead}),
we may choose these $x_i$, $i=1,2,3$, to
lie in $q^{-1}_0(R'(F_3,H))$.

Consider the orbit
$\text{Orb}_{_{K^3}}(x_1\, , x_2\, , x_3)$, of $(x_1\, , x_2\, ,
x_3)\,\in\, K^3$ for the adjoint action of the group $K^3$ on itself.
Clearly we have the following identification of this orbit:
\[
\text{Orb}_{_{K^3}}(x_1\, , x_2\, , x_3) \,= \,
(K/{C(x_1)}) \times (K/{C(x_2)}) \times (K/{C(x_3)})
\]
with $C({x_i})\, \subset\, K$ being the centralizer of $x_i$.
The conjugation action of $K$ on $\text{Hom}(F_3,K)$ coincides
with the restriction of above action of $K^3$ to the image
of the diagonal map $K \, \hookrightarrow\, K^3$. Therefore,
the fibre $K/Z$ through the point
$(x_1, x_2, x_3)\,\in\, q^{-1}_0(R'(F_3,H))$ of the $K/Z$-bundle
$$
q^{-1}_0( R'(F_3,H)) \,\longrightarrow\, R'(F_3,H)
$$
is contained in the orbit
\[
(K/{C(x_1)}) \times (K/{C(x_2)}) \times (K/{C(x_3)})\, .
\]

Now if we choose a point $(x_1, x_2, x_3) \in K^3$
which is general enough, then by the definition of the inverse image
$q^{-1}_0( R'(F_3,H))$ and Lemma \ref{newstead}(2)
it follows immediately that the complement of the open dense subset
\[
q^{-1}_0(R'(F_3,H))\bigcap (\frac{K}{C(x_1)} \times
\frac{K}{C(x_2)} \times \frac{K}{C(x_3)})\,\subset\,
\frac{K}{C(x_1)} \times\frac{K}{C(x_2)}\times\frac{K}{C(x_3)}
\]
is of codimension at least four.

Since the image of $K/Z$ in $q^{-1}_0( R'(F_3,H))$ lies in $(q^{-1}(
R'(F_3,H)))\bigcap \text{Orb}_{_{K^3}}(x_1, x_2, x_3)$, whose
complement is of codimension at least four in
\[
\text{Orb}_{_{K^3}}(x_1, x_2, x_3)
\,=\, (K/{C(x_1)}) \times (K/{C(x_2)}) \times (K/{C(x_3)})
\]
it follows
that the homomorphism $ \pi_1(K/Z) \,\longrightarrow\, \pi_1(q^{-1}_0
(R'(F_3,H)))$ in the proposition factors through
\[
\pi_1((q^{-1}_0( R'(F_3,H))) \cap \text{Orb}_{_{K^3}}(x_1, x_2, x_3))
\]
\[
=\,
\pi_1((K/{C(x_1)}) \times (K/{C(x_2)}) \times (K/{C(x_3)})) \,=
\,\pi_1(K/T)^3
\]
where $T$ is a maximal torus of $K$.

For any maximal torus $T\, \subset\, K$, the quotient $K/T$ is
diffeomorphic to $H/B$, where $B$ is a Borel subgroup of $H$.
Since $H/B$ is simply connected, we conclude that
$(K/{C(x_1)}) \times (K/{C(x_2)}) \times (K/{C(x_3)})$ is simply
connected. This completes the proof of the proposition.
\end{proof}

{\it We now complete the proof of the non-existence of the covering
${\widetilde M} \longrightarrow q^{-1}( R'(F_3,H))$
as in \eqref{eqdig2}.}

Since the homomorphism of fundamental groups
$ \pi_1(K/Z) \,\longrightarrow\, \pi_1(q^{-1}_0( R'(F_3,H)))$ is
trivial (see Proposition \ref{claim}),
the induced covering $K\,\longrightarrow\, K/Z$ is
trivial (see the diagram \eqref{eqdig2}).
But this is a contradiction to the fact that $K$ is connected.

Therefore, we have proved the following theorem:

\begin{theorem}\label{nonexist}
If the centre $Z(H)$ is non-trivial, then there is no universal
bundle over $X\times {\mathcal M}'_X(H)$.
\end{theorem}


\end{document}